\documentclass[12 pt]{amsart}
\usepackage{amssymb}
\usepackage{enumerate}
\usepackage{mathrsfs}
\usepackage[usenames,dvipsnames,table,xcdraw]{xcolor}
\usepackage{graphicx}
\usepackage{multicol}
\usepackage{todonotes}
\usepackage{etoolbox}
\usepackage{wrapfig}
\usepackage{enumitem}
\usepackage{sidenotes}
\usepackage{accents}

\makeatletter
\newcommand*{\dotprod}{}
\DeclareRobustCommand*{\dotprod}{%
  \mathbin{\mathpalette\dotprod@{}}%
}
\newcommand*{\dotprod@scalefactor}{.6}
\newcommand*{\dotprod@widthfactor}{1.15}
\newcommand*{\dotprod@}[2]{%
  \sbox0{$#1\vcenter{}$}
  \sbox2{$#1\cdot\m@th$}%
  \hbox to \dotprod@widthfactor\wd2{%
\hfil
\raise\ht0\hbox{%
  \scalebox{\dotprod@scalefactor}{%
\lower\ht0\hbox{$#1\bullet\m@th$}%
  }%
}%
\hfil
  }%
}

\theoremstyle{plain}

\theoremstyle{definition}

\author{Benjamin Richeson}
\address{SUNY Polytechnic Institute}
\email{benricheson101@gmail.com}

\author{David Richeson}
\address{Department of Mathematics and Computer Science, Dickinson College}
\email{richesod@dickinson.edu}

\title{What's the Best Seat in the Game Left, Center, Right?}
\date{\today}   
                                    
\begin{document}
\begin{abstract}
\textit{Left, Center, Right} is a popular dice game. We analyze the game using Markov chain and Monte Carlo methods. We compute the expected game length for two to eight players and determine the probability of winning for each player in the game. We discuss the surprising conclusions of which players have the highest and lowest chance of winning, and we propose a small rule change that makes the game a little more fair.
\end{abstract}
\maketitle

\textit{Left, Center, Right} is a popular dice game typically played by at least three players. Its gameplay is purely probabilistic, giving the players no opportunities to make any decisions. 

The game begins with the players sitting in a circle with three chips each (or, for those more adventurous, three \$1, \$5, or \$10 bills!). There are three special six-sided dice. Three faces of each die have dots, one face has a star, one has the word LEFT (henceforth just L), and one has the word RIGHT (R). For simplicity, we will refer to a dot as HOLD (H) and the star as CENTER (C). To play this game with ordinary dice, simply assign the numbers 1 through 6 to the four possibilities, such as treating the even numbers as HOLDs, 1 as LEFT, 3 as RIGHT, and 5 as CENTER.

On a player's turn, they roll as many dice as they have chips (up to a maximum of three) and act based on the dice rolls. For each die that comes up L, R, or C, the player passes one chip to the person on their left, one to the person on their right, or one into the pot in the center of the table, respectively. Dice that come up H require no action.

The play continues clockwise in this fashion. If a player has no chips, they do not roll, but they may not be out of the game because someone may pass them a chip later. The chips in the center stay in the center. The last player to have any chips wins the game (or, if playing with money, wins the pot of cash). 

Because dice rolls completely determine the game, it lends itself to study using the theory of Markov chains. A Markov chain is a stochastic model in which, at a given state, the probability of reaching another state is fixed and does not depend on previous events. \textit{Chutes and Ladders}---in which gameplay is determined by the fixed locations of chutes and ladders on the board and spins of the spinner---is another familiar game that can be modeled using Markov chains. For a variety of examples, see \cite{Abbott:1997,Althoen:1993,Ash:1972,Connors:2014,Crans:2015,Diaconis:2001,Do:2025,Lakshtanov:2012,Osborne:2003,Pierce:2015}. For an analysis of a variant of \textit{Left, Center, Right} that does not use Markov chains, see \cite{Torrence:2016}. No previous exposure to Markov chains is required to follow the arguments in this article, but to learn more about this area of probability, see, for instance, \cite{Roberts:1976}.

Markov chains allow us to obtain precise answers to questions like: What is the expected length of an $n$-player game of \textit{Left, Center, Right} (in terms of the number of dice rolls)? What is the probability each player will win? What is each player's expected number of chips after $m$ turns? And so on.

Computing these answers using Markov chains is elegant and exact, but the more players, the more computationally expensive. Another option is to compute approximations of these values using Monte Carlo methods: Simply play the game many times and report the approximate values. 

In this article, we investigate \textit{Left, Center, Right} using these techniques. Ultimately, our question is, where should you sit to maximize your chance of winning the game? 

\section{An analysis using Markov chains}
\label{sec:Markov}
To study \textit{Left, Center, Right} mathematically, we first define the state of the game. If there are $n$ players, we can represent the game state by an $(n+1)$-tuple of non-negative integer values, $(c_1,c_2,\ldots,c_n;j)$, where $c_i$ is the number of chips in front of player~$i$, and $j$ indicates whose turn it is. 

\begin{figure}[ht]
\begin{center}
\includegraphics{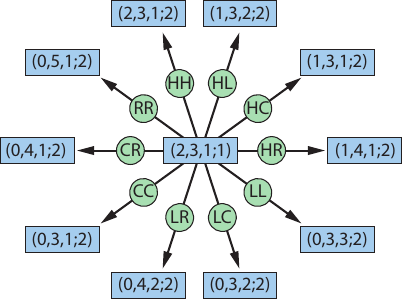}
\end{center}
\caption{When the players~1 through 3 have two, three, and one chips, respectively, and it is player~1's turn, there are ten possible outcomes.}
\label{fig:Markov}
\end{figure}

If the game is in a particular state, there are well-defined probabilities for transitions to other states. For instance, if the game is in state $(2,3,1;1)$, then player~1 has two chips, player~2 has three, player~3 has one, and it is player~1's turn to roll. Since player~1 has two chips, they must roll two dice. There are ten possible rolls: HH, HL, HC, HR, LL, LC, LR, CC, CR, and RR. Figure~\ref{fig:Markov} shows the ten possible states that follow $(2,3,1;1)$. Notice that with each L, player~1 passes a chip to player~3; with each R, they pass to player~2, each C decreases their chips by one, and each H has no effect. Also, moving to the next state increments the player counter by 1 mod 3.

\begin{figure}[ht]
\begin{center}
\includegraphics{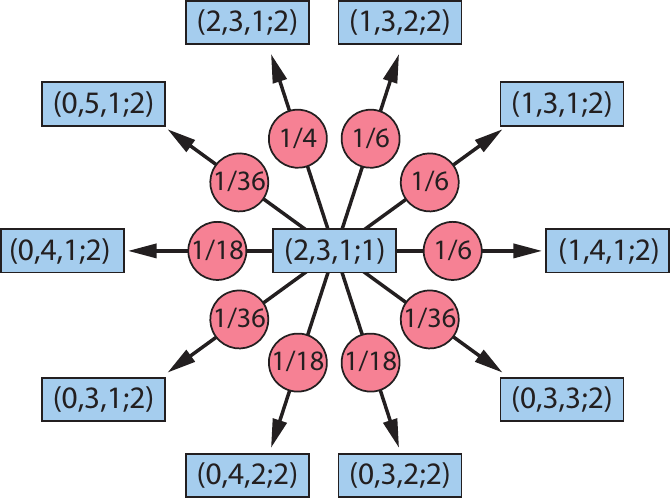}
\end{center}
\caption{These are the transition probabilities when the players~1 through 3 have two, three, and one chips, respectively, and it is player~1's turn.}
\label{fig:Markovprobs}
\end{figure}

Figure~\ref{fig:Markovprobs} shows the probabilities associated with each transition. We will say more about these probabilities shortly.

Once we know the transition probabilities between all the states, we obtain a directed graph with probabilities on each edge. Playing the game corresponds to beginning at the state $(3,3,\ldots,3;1)$ and traversing the graph, making decisions based on the probabilities on the edges. The game ends when we get to a state $(c_1,c_2,\ldots, c_n;j)$ in which only one of the $c_i$ is nonzero. In the directed graph, such vertices have a loop from the vertex to itself labeled with a probability of 1.

We can say exactly how many states are in the graph for an $n$-player game. There are $3n$ chips that must be divided into $n+1$ piles---the $n$ players and the center region.  This problem is equivalent to the well-known ``stars and bars'' question---how many ways can $3n$ stars (*) and $n$ bars $(|)$ be arranged? For instance, if $n=4$, we may have the following arrangement.
\[*|****|***|~|****\]
This corresponds to player~1 having one chip, player~2 having four, player~3 having three, player~4 having none, and the center having four. 

We must choose locations for the $n$  bars out of the $3n+n$ possible locations; the number of such possibilities is \[{3n+n\choose n}={4n\choose n}.\] 
However, we must subtract 1 from this quantity because no game state has every chip in the center. Also, in addition to keeping track of the distribution of chips, we must keep track of whose turn it is. Thus, the total number of states is 
\[n\left({4n\choose n}-1\right)=n\cdot{4n\choose n}-n.\] 
We apply this equation to the first seven $n$-values to obtain the values in Table~\ref{tab:numvertices}.

\begin{table}[hbt]
\begin{center}
\begin{tabular}{cc}
\hline
Number of players ($n$)& Number of states \\
\hline
2& 54\\
3& 657 \\
4& 7,276\\
5& 77,515 \\ 
6& 807,570\\
7 & 8,288,273\\
8 & 84,146,392\\
\hline
\end{tabular}
\end{center}
\caption{The number of vertices of the graph for the game with $n$ players.}
\label{tab:numvertices}
\end{table}

Notice that the number of states increases quickly---by about a factor of 10 with each new player. In fact, for large $n$, we have  
\begin{align*}
\frac{\text{number of states for }n+1\text{ players}}{\text{number of states for }n\text{ players}}
&=\frac{(n+1)\cdot{4(n+1)\choose n+1}-(n+1)}{n\cdot{4n\choose n}-n}\\
&\approx 
\frac{(n+1)\cdot{4(n+1)\choose n+1}}{n\cdot{4n\choose n}}\\
&=\frac{(4n+4)(4n+3)(4n+2)(4n+1)}{n(3n+3)(3n+2)(3n+1)}\\
&\approx \frac{256}{27}\approx 9.5.
\end{align*}

Next, we compute the probabilities on each edge. Suppose we want to compute the probability of rolling three dice and getting two R's and one H. For one die, the probability of rolling R is $\frac16$, and the probability of rolling H is $\frac12$. But the order of rolling them is unimportant, and the two R's are identical. So, the probability is
\[\frac{3!}{2!1!}\cdot\left(\frac16\cdot\frac16\cdot\frac12\right)=\frac{1}{24}.\] Similarly, the probability of rolling one L, one C, and one H is \[\frac{3!}{1!1!1!}\cdot\left(\frac16\cdot\frac16\cdot\frac12\right)=\frac{1}{12}.\] The other probabilities are computed similarly and are shown in Table~\ref{tab:probabilities}.

\begin{table}[hbt]
\bgroup
\def\arraystretch{1.5}
\begin{center}
\begin{tabular}{cccc}
\hline
Roll Outcomes&&&\\
(Omitting H) & One Die & Two Dice & Three Dice  \\
\hline
No L, C, or R rolls & $\frac12$ & $\frac12\cdot\frac12=\frac14$ &$\frac12\cdot\frac12\cdot\frac12=\frac18$\\
L, C, or R& $\frac16$ &$2(\frac16\cdot\frac12)=\frac16$ & $3(\frac12\cdot\frac12\cdot\frac16)=\frac18$\\
LC, LR, or CR& --- & $2(\frac16\cdot\frac16)=\frac{1}{18}$ &$6(\frac16\cdot\frac16\cdot\frac12)=\frac{1}{12}$\\
LL, RR, or CC&--- &$\frac16\cdot\frac16=\frac{1}{36}$ &$3(\frac16\cdot\frac16\cdot\frac12)=\frac{1}{24}$\\
LCR&--- &--- &$6(\frac16\cdot\frac16\cdot\frac16)=\frac{1}{36}$\\
\begin{tabular}{@{}c@{}}LLC, LLR, LRR, \\ CRR, LCC, or CCR\end{tabular} & ---&--- &$3(\frac16\cdot\frac16\cdot\frac16)=\frac{1}{72}$\\
LLL, RRR, CCC & ---&--- & $\frac16\cdot\frac16\cdot\frac16=\frac{1}{216}$\\
\hline
Number of Outcomes&4&10&20\\
\hline
\end{tabular}
\end{center}
\egroup
\caption{The probabilities of the various roll outcomes when throwing one, two, or three dice and the number of possibilities.}
\label{tab:probabilities}
\end{table}

Note also that if it is player~$i$'s turn but they have no chips, then with probability 1, the game transitions to the next state; chip counts remain the same, and it is player~$(i+1)$'s turn to roll.

The directed graph labeled with probabilities gives us a nice visual understanding of the Markov process. But mathematically, collecting this information in matrix form is more useful. The \emph{transition matrix} associated to a Markov chain with $m$ states is an $m\times m$ matrix in which the $(i,j)$-entry is the probability of going from the $i$th state to the $j$th state. Notice that in such a matrix, all entries are nonnegative, and every row sums to 1.

In our graph, it is possible to find a path from any vertex to a winning state. Such a Markov chain, in which there is a path from any vertex to some state in which the only outgoing edge is a loop, is called an \emph{absorbing Markov chain}, and the vertices whose only outgoing edge is a loop are \emph{absorbing states}. Nonabsorbing states are called \emph{transient states.} Vertex $i$ is absorbing if and only if the $i$th row in the transition matrix consists of all zeros except a 1 in column $i$. 

For simplicity, we number the vertices so that the transient states are listed first and the absorbing states are listed last. Let's say that the Markov chain has $t$ transient states and $a$ absorbing states. Then, the transition matrix has the form 
\[
P=\left(\begin{array}{c|c}
Q  & R\\\hline
 \mathbf{0} & I_a 
\end{array}\right)
,\]
where $Q$ is a $t\times t$ matrix, $R$ is a $t\times a$ matrix, $\mathbf{0}$ is an $a\times t$ matrix of 0's, and $I_a$ is the $a\times a$ identity matrix.

The \emph{fundamental matrix} for the Markov chain is \[N=\left(I_t-Q\right)^{-1},\] in which $I_t$ is the $t\times t$ identity matrix. Remarkably, the sum of the values in the $i$th row of $N$ is the expected number of steps before being absorbed when starting in transient state $i$. Thus, if vertex 1 is the starting vertex, corresponding to state $(3,3,\ldots,3;1)$, the sum of the entries in the first row is the expected number of dice rolls until the end of the game.

We can even use the fundamental matrix to compute the standard deviation of this value. Let \[\mathbf{w}=N\left(\begin{smallmatrix}
1  \\
1 \\
\vdots \\
1
\end{smallmatrix}
\right)=\left(\begin{smallmatrix}
w_1  \\
w_2 \\
\vdots \\
w_t
\end{smallmatrix}
\right)\text{ and }
\mathbf{w}_{\text{sq}}=\left(\begin{smallmatrix}
w_1^2  \\
w_2^2 \\
\vdots \\
w_t^2
\end{smallmatrix}
\right).
\]
Then, the variance for the number of steps before being absorbed, when starting in transient state $i$, is the $i$th entry of \[\left(2N-I_t\right)\mathbf{w}-\mathbf{w}_{\text{sq}}.\] Hence, the standard deviation is the square root of that value.

For instance, we computed that the average game length for three players is 18.9 rolls with a standard deviation of 8.1 rolls. 

The fundamental matrix is also useful for analyzing Markov chains with multiple absorbing states, such as ours. In this case, the $(i,j)$-entry of the matrix \[B=NR\] is the probability of the game ending in absorbing state $j$ when starting from transient state $i$. Assuming, again, that we begin at vertex 1, then the sum of all the entries $(1,j)$ of $B$ in which $j$ corresponds to a victory for player~$p$ is the probability that player~$p$ wins.  

For the three-player version of \textit{Left, Center, Right}, we found that the probabilities of victory for players 1, 2, and 3 are 0.307, 0.328, and 0.365, respectively. Thus, to ensure the best chance of winning, you should choose to be player~3. 

In the Results section, we give the expected game lengths and the probabilities of victory for other numbers of players.

\section{Monte Carlo simulation}
\label{sec:MonteCarlo}
Because the number of states increases rapidly with increasing $n$, computing the desired quantities using Markov chains is impractical or impossible for large $n$. In these cases, we can use the Monte Carlo method. We have the computer toss the dice and repeatedly play the game. Then, using this data, we compute the average game length and the percentage of the games that each player wins. 

We know that the more games we play, the better our approximation, but the question is, how accurate are our approximations if we play the games some specific number of times? What if we played 1,000 times? One million times? One billion times?

Suppose we play $R$ games in our Monte Carlo simulation. When we do, we obtain an approximation $\overline{x}$ for the average number of turns to win or the fraction of the games won by one of the players. By the law of large numbers, the larger $R,$ the closer $\overline{x}$ is to the desired quantity. Moreover, by the central limit theorem, the results of such Monte Carlo simulations follow a roughly normal distribution centered around the desired value. In this scenario, the \textit{Monte Carlo error}, $\sigma_R$, is the standard deviation of the Monte Carlo estimator taken across hypothetical repetitions of the simulation, where each simulation consists of $R$ replications. (See \cite{Koehler:2009} or \cite[\S 3.4]{Suess:2010}, for instance.) If $\hat{\sigma}$ is the calculated standard deviation obtained from the $R$ games, \[\sigma_R\approx\frac{\hat{\sigma}}{\sqrt{R}}=\frac{1}{\sqrt{R}}\sqrt{\frac{1}{R-1}\sum_{i=1}^R(x_i-\overline{x})^2}\approx\sqrt{\frac{1}{R^2}\sum_{i=1}^R(x_i-\overline{x})^2}.\]
Then, for instance, the absolute error of our approximation $\overline{x}$ is larger than two standard deviations, $2\sigma_R$, in only five percent of simulations.

We played $R=100,000,000$ five-player games of \textit{Left, Center, Right}. The average game length was $\overline{x}=49.93574$, and the standard deviation was $\hat{\sigma}=15.94761$. Thus, \[2\sigma_R\approx\frac{2\hat{\sigma}}{\sqrt{10^8}}\approx0.00319.\]
So, the true average game length is very likely between 49.9326 and 49.9389.

Obtaining multiple decimal places of accuracy for the expected game length is unnecessary, but we do require a high degree of accuracy when trying to see which players are most and least likely to win. While we could use the Monte Carlo error formula to compute the error, there is a more straightforward method. 

Suppose we've played $R$ simulated games, and player~$k$ wins $w$ of them. Then $\hat p=w/R$ is the sample proportion of games won. If $p$ is the true proportion, the standard error is 
\begin{align*}
\sigma_{\hat p}
&=\sqrt{\frac{p(1-p)}{R}}\\
&\approx\sqrt{\frac{\hat p(1-\hat p)}{R}}\\
&=\sqrt{\frac{1}{R}\left(\frac{w}{R}\right)\left(1-\frac{w}{R}\right)}\\
&=\sqrt{\frac{Rw-w^2}{R^3}}=\hat \sigma_{\hat p}.
\end{align*}
Moreover, the Central Limit Theorem tells us that $(\hat p-p)/\hat\sigma_{\hat p}$ converges in distribution to a standard normal random variable.

For instance, in our Monte Carlo simulation of five-player games, player 1 won $w=19,375,770$ of the $R=100,000,000$ games. So $\hat{p}= 0.1937577.$ Then,
\[2\hat \sigma_{\hat p}=2\sqrt{\frac{(10^8)(19375770)-(19375770)^2}{(10^8)^3}}\approx0.000079048\]
So, the probability that player~1 wins is very likely between 0.19368 and 0.19384.

\section{Results}
\label{sec:Results}

The game \textit{Left, Center, Right} comes with 24 chips, allowing a maximum of eight players, although there is no upper limit for the number of players (other than time and patience). The game is recommended for three players or more, but two could play the game as well (a player's only opponent would be on both the right and left). So, using the Markov chain techniques described earlier, we computed the expected game length and the probability that each player wins for an $n$-player game for $n=2, 3, 4$ players. Then, using the Monte Carlo method, we computed the same quantities for $n=5,6,7,8$ players (using $R=100,000,000$). These values are collected in Tables~\ref{tab:numturns} and \ref{tab:probwin}.

\begin{table}[hbt]
\begin{center}
\begin{tabular}{ccc}
\hline
Number of players ($n$) &Expected number of turns & Standard deviation \\
\hline
2 & 5.8 & 3.7 \\ 
3 & 18.9 & 8.1\\
4 & 33.9& 12.1\\
5 & 49.9& 15.9 \\
6 & 66.6 & 19.6\\
7 & 84.0& 23.3\\
8 & 101.9 & 27.1 \\
\hline
\end{tabular}
\end{center}
\caption{The expected number of turns in an $n$-player game. }
\label{tab:numturns}
\end{table}

\begin{table}[hbt]
\begin{center}
\begin{tabular}{ccccccccc}
\hline 
\multicolumn{1}{c}{\begin{tabular}{@{}c@{}}Number of\end{tabular}} & \multicolumn{8}{c}{Probability that player~$k$ wins}   \\
players ($n$) & 1 & 2 & 3 &4 & 5 & 6 & 7 & 8   \\
\hline
2& \cellcolor{red!40}0.382& \cellcolor{ForestGreen!40}0.618 & & & & & & \\
3& \cellcolor{red!40}0.307 & 0.328 & \cellcolor{ForestGreen!40}0.365 & & & & & \\
4& \cellcolor{red!40}0.239 & 0.243 & 0.255 & \cellcolor{ForestGreen!40}0.262 & & & & \\
5& 0.194 & \cellcolor{red!40}0.194 & 0.200 & \cellcolor{ForestGreen!40}0.206 & 0.206 & & & \\
6& 0.162 & \cellcolor{red!40}0.161 & 0.164 & 0.169 & \cellcolor{ForestGreen!40}0.173 & 0.170 & & \\
7& 0.139 & \cellcolor{red!40}0.139 & 0.140 & 0.143 & 0.147& \cellcolor{ForestGreen!40}0.148 & 0.145 & \\
8& 0.122 & \cellcolor{red!40}0.120 & 0.121 & 0.124 & 0.127 & 0.129 & \cellcolor{ForestGreen!40}0.130 & 0.127 \\
\hline
\end{tabular}
\end{center}
\caption{The probability that player~$k$ wins in an $n$-player game. The red cells indicate the players with the lowest probability and the green cells are the players with the highest probability.}
\label{tab:probwin}
\end{table}

The results in Table~\ref{tab:probwin} surprised us. We expected that the chance of winning would increase the later the player was in the lineup because each time a player rolls, they will most likely give away chips. But that pattern did not hold. For the larger $n$-values, player~1 did better than player~2, and player~$n$ did worse than player~$n-1.$ 

\begin{figure}[ht]
\includegraphics{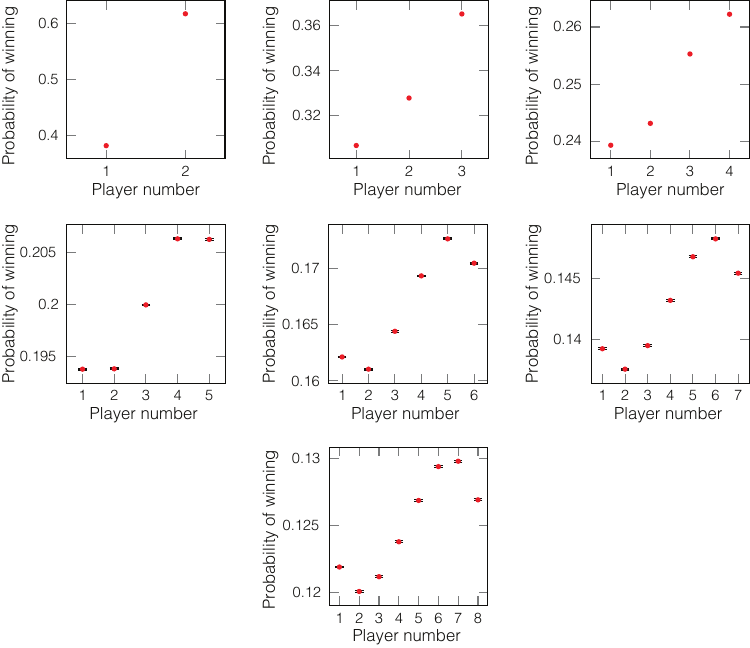}
\caption{Probability of victory for games of two to eight players (with error bars for those computed using Monte Carlo methods).}
\label{fig:SevenGraphs}
\end{figure}

\begin{figure}[ht]
\begin{center}
\includegraphics{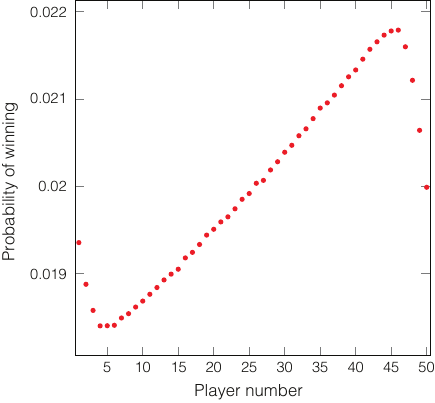}
\end{center}
\caption{The probability of victory for 50 players.}
\label{fig:fifty}
\end{figure}

Figure~\ref{fig:SevenGraphs} shows graphs of the probabilities with error bars, which are barely distinguishable from the points themselves. We see that the effect is real. In fact, we simulated 100 million games with 50 players and saw that the same behavior occurred for the first and last several players (see Figure~\ref{fig:fifty}).

\section{Discussions and conclusions}
Before we discuss the surprising behavior exhibited by the data, we will make some general observations about the game. First, since one player cannot take chips from another player, a player cannot win the game on their own turn. Said another way, a player wins the game when the only remaining player loses by putting their remaining chips in the center. And, in fact, a player does not want to roll. The best outcome for a player on their turn is that they roll only H's and end with the same number of chips they began with; more often, they give away chips.

When a player does roll, we can say, on average, how many chips they will lose and where they will go. When rolling one die, we expect a player to pass $1/6$ of a chip to the left, $1/6$ to the right, and $1/6$ to the center. When rolling two dice, they can expect to pass $1/3$ to the left, $1/3$ to the right, and $1/3$ to the center. When rolling three, they expect to pass $1/2$ to the left, $1/2$ to the right, and $1/2$ to the center. One consequence of this observation is that it is better to have rich neighbors than poor neighbors---probabilistically speaking, the more dice your neighbor rolls, the better it is for you.

The game begins in an unusual fashion. The first player will have three chips on their first turn, but the rest will, on average, have more than three on their first turn since the neighbor on their right may have passed them some. Moreover, the last player, player~$n$, may receive chips from two players before they roll---from players~1 and $n-1$. For more consistency, one could make one simple rule change: When player 1 rolls for the first time, they treat any L's as H's. That is, they hold the chips they would otherwise pass to player $n$. We will call this the modified start rule. 

\begin{figure}[ht]
\begin{center}
\includegraphics{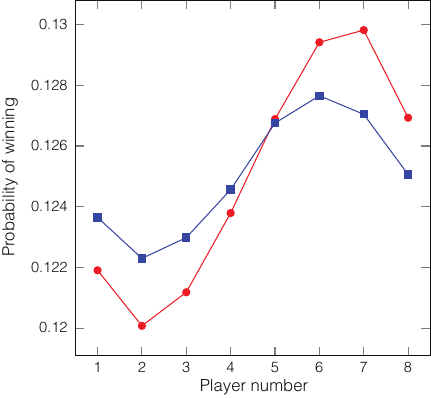}
\end{center}
\caption{Probability of victory in an eight-player game under the ordinary (red circles) and modified (blue squares) rule.}
\label{fig:modifiedrules}
\end{figure}

Although the modified start rule does not make the game fair, it does make it more fair. As we see in Figure~\ref{fig:modifiedrules}, for eight players, the difference between the largest and smallest winning probabilities is half the size with this new rule. Note that the probabilities still have the same S-shape.

In the following discussion, we will assume the modified start rule since it is slightly easier to see how the graphs get their shapes. Also, for concreteness, we will assume there are six players.

Table~\ref{tab:firstround} shows each player's expected number of chips on the first 10 turns. The numbers in the first seven rows are exact since the players will roll three dice and will, on average, give half a chip to the player on their left, half a chip to the player on their right, and half a chip to the center. We computed the rest of the values using Monte Carlo methods.

\begin{table}[ht]
\begin{center}
\begin{tabular}{cccccccc}
\hline
Turn&\multicolumn{6}{c}{Player ($k$)}\\
($j$)& 1 & 2& 3 & 4 & 5 & 6\\
\hline
1&\cellcolor{blue!40}3&3&3&3&3 & 3\\
2&2$^*$&\cellcolor{blue!40}3.5&3&3&3 &3$^*$ \\
3&2.5&2&\cellcolor{blue!40}3.5&3&3 &3\\
4&2.5&2.5&2&\cellcolor{blue!40}3.5&3 &3\\
5&2.5&2.5&2.5&2&\cellcolor{blue!40}3.5 &3\\
6&2.5&2.5&2.5&2.5&2 &\cellcolor{blue!40}3.5 \\
7&\cellcolor{blue!40}3&2.5&2.5&2.5&2.5 &2\\
8&1.77& \cellcolor{blue!40}2.91&2.5&2.5&2.5&2.41\\
9&2.18&1.70& \cellcolor{blue!40}2.90&2.5 &2.5&2.41\\
10&2.18&2.10&1.69& \cellcolor{blue!40}2.90&2.5&2.41\\
\end{tabular}
\end{center}
\caption{The expected number of chips for each player at the start of turn $j$ ($^*$assuming the modified start rule). The blue cell indicates which player will roll next.}
\label{tab:firstround}
\end{table}

Two factors produce the S-shape. First, the only person able to end the game is the person rolling (by losing the game). In that respect, the later a player is in the ordering, the better. If players don't lose on their turn, they are safe until their next turn. This gives the graphs a general upward trend. 

Looking closer at the expected values in Table~\ref{tab:firstround}, we see a ``wave'' traveling around the table. In the first round, for instance, after player~$k$ rolls, player~$k+1$ is the ``crest'' of the wave since they have, on average, 3.5 chips; player~$k$ is the ``trough'' of the wave with two chips; and player~$k-1$ reaches 2.5 where it remains for the rest of the round. However, something slightly different happens in the transition between rounds. When it is player~1's second turn, they have, on average, three chips. But that's just the expected value, so they will not always roll three dice. On average, they pass 0.41 chips to the left, center, and right. Player~6 ends up with an average of 2.41 chips, not 2.5 like the other players in that round. Moreover, after player~2 rolls, player~1 ends up with an expected number of 2.18 chips, which is larger than the expected values of the later players in this round (2.10 chips). Thus, we see how this transition from one round to the next transfers chips from the last player to the first player and produces the S-shape.

If there is one strategic takeaway from this analysis, it is that when playing \textit{Left, Center, Right}, try to be one of the last players to roll rather than one of the first players to roll. Table~\ref{tab:probwin} gives precise seating advice in games of two to eight players. The difference between the best and worst seats in a three-player game is quite dramatic, producing a best-case 36.5\% chance of winning versus a worst-case 30.7\% chance. As the number of players increases, the margins decrease. The best seat in an eight-player game has a 13\% chance of winning, which is only slightly larger than the worst seat, which has a 12\% chance. 

Still, if each of the eight players starts with \$30, player~7's expected winnings of $\$30\cdot 8\cdot 0.13-\$30=\$1.20$ is quite a bit better than player~2's expected loss of the same amount.\\

\noindent\textbf{Acknowledgment.} We thank Jeff Forrester for his helpful comments about Monte Carlo error. 

\bibliographystyle{plain} 
\bibliography{LCR.bib}

\end{document}